\documentstyle[12pt]{article}

\def\be{\begin{equation}}
\def\ee{\end{equation}}
\def\bea{\begin{eqnarray}}
\def\eea{\end{eqnarray}}

\def\1{\'{\i}}

\begin{document}

\thispagestyle{empty}

 \
\hfill \

\ 
\vspace{2cm}

\begin{center}
{\Large{\bf{Discretizations of the Schr\"odinger equation}}}

{\Large{\bf{with quantum algebra symmetry}}} 
\end{center}

\bigskip\bigskip

\begin{center} Angel Ballesteros, Francisco J.
Herranz\footnote{Communication presented in the 5th Wigner Symposium,
Vienna, 1997}  and Preeti Parashar 
\end{center}

\begin{center} {\it { Departamento de F\1sica, Universidad
de Burgos} \\   Pza. Misael Ba\~nuelos, 
E-09001-Burgos, Spain}
\end{center}

\bigskip\bigskip\bigskip

\begin{abstract}
Two quantum Hopf structures for the Schr\"odinger algebra as well  as
their corresponding differential-difference realizations are presented.
For each case a (space or time) discretization of the Schr\"odinger
equation is deduced and the quantum Schr\"odinger generators are shown
to be symmetry operators.
\end{abstract}

\bigskip\bigskip\bigskip

\newpage

\section{A   space discretization of the Schr\"odinger
equation}

The (1+1)-dimensional Schr\"odinger  algebra \cite{Hagen,Nied} ${\cal S}$
is generated by the time translation $H$, space translation $P$, Galilean
boost $K$, dilation $D$, conformal transformation $C$  and mass
operator $M$. The  non-zero Lie brackets of ${\cal S}$ are:
\be
\begin{array}{llll}
 [D,P]=-P \quad  &[D,K]=K \quad &[K,P]=M \quad &[K,H]=P \cr
 [D,H]=-2H \quad  &[D,C]=2C  \quad &[H,C]=D \quad  &[P,C]=-K  . 
\end{array}
\label{aa}
\ee
The usual differential realization
of the Schr\"odinger generators in terms of the time and space
coordinates $(t,x)$ and two constants $a$, $m$, is given by:
\be
\begin{array}{l}
 H=\partial_t \qquad  P=\partial_x \qquad M=m \qquad
 K=-t \partial_x - m x\cr 
D=2 t \partial_t +
x\partial_x -a \qquad  C=t^2 \partial_t + t x \partial_x - a t  + \frac m2
x^2 .
\end{array}
\label{ab}
\ee
The heat/Schr\"odinger equation (SE) can be  
obtained from this representation  as the Casimir  
$E=P^2 -2 M H$ of the Galilei subalgebra $\{K,H,P,M\}$ of ${\cal S}$:
\be
(\partial_x^2 - 2 m \partial_t) \phi(x,t)=0.
\ee 
We   say that an operator $S$ is a symmetry of $E$ if $S$
transforms solutions into solutions:
$(ES) \phi(x,t)= (\Lambda S)\phi(x,t)$, 
where   $\Lambda$ is
another operator. The  Galilei generators are
symmetries $S$ of the SE  with $\Lambda=S$. The dilation
is also a symmetry, since (\ref{aa}) implies that $[E,D]=2
E$. Finally, we compute $[E,C]=-(K P + P K +2 M D)$ and by introducing
the representation  (\ref{ab})
with  $a=-1/2$  we obtain that $[E,C]=2 t E$.

Let $U_z^{(s)}({\cal S})$ be the quantum  Schr\"odinger algebra
\cite{twoa} whose (nonstandard) classical $r$-matrix is
$r=z P\wedge (D +\frac 12 M)$; its
 coproduct and  commutators  are
\be
\begin{array}{l}
 \Delta(P)=1\otimes P + P \otimes 1 \qquad
\Delta(M)=1\otimes M + M \otimes 1 \cr
 \Delta(H)=1\otimes H + H\otimes e^{-2zP}\cr
 \Delta(K)=1\otimes K + K\otimes e^{zP}-z D'\otimes
e^{zP}M\cr
 \Delta(D)=1\otimes D + D\otimes e^{zP} + \frac 12  M\otimes
(e^{zP}-1) \cr
 \Delta(C)=1\otimes C + C\otimes  e^{2zP}+\frac z2 K\otimes  e^{zP} D' 
-\frac z2  D'\otimes  e^{zP}(K+z D' M) ,
\end{array} 
\label{ba}
\ee
\be
\begin{array}{l}
 [D,P]=\frac 1z ({1- e^{z P}}) \quad\    [D,K]=K\quad\ 
 [K,P]=M  e^{z P } \quad\  [M,\,\cdot\,]=0 \cr
  [D,H]= - 2 H \qquad   [D,C]=2 C -\frac z2 K D' 
\qquad  [H,P]=0 \cr
 [H,C]= \frac 12(1+  e^{-z P}) D' -\frac 12 M +  z K H  \qquad 
[K,C]=-\frac z2 K^2 \cr
  [P,C]=-\frac 12 (1+ e^{zP}) K - \frac z2   e^{zP} M D'
  \qquad   [K,H]=\frac 1z({1-e^{-zP}})  , 
\end{array} 
\label{bb}
\ee
where $D'=D+\frac 12 M$.
We   deduce the   quantum SE with $U_z^{(s)}({\cal S})$ as its
quantum symmetry  algebra.  A differential-difference realization of
(\ref{bb}) reads
\be
\begin{array}{l}
 P=\partial_x \qquad H=\partial_t  \qquad M= m \cr
 K=- t \frac 1z ({1- e^{-z \partial_x}})   - m x e^{z \partial_x}
\qquad  
 D=2  t  \partial_t  +
x \frac 1z({  e^{z \partial_x}-1}) -a   \cr
 C= t^2 \partial_t  e^{-z \partial_x} + t x \left( \frac 1z {\sinh
 z\partial_x}  - z m \partial_t  e^{z \partial_x}\right) +\frac 12 m x^2
e^{z
\partial_x}\cr
 \qquad -\frac 12 t \left\{ a (1+e^{-z \partial_x}) + 
(1+ \frac 12 m)(1-e^{-z \partial_x})\right\} \cr
 \qquad 
  +\frac 12 z m (1+a- \frac 12 m) x 
e^{z \partial_x} .
\end{array} 
\label{bc}
\ee
At the level of commutation rules (\ref{bb}), the Galilei structure remains
as a subalgebra and its deformed Casimir   is
$ E_z = \frac 1{z^2} ({1 -  e^{- z P}})^2 - 2 M H $.
We define a space discretization of the SE as the action of $E_z$ on
$\phi(x,t)$ through (\ref{bc}), that is, 
\be
 \left(    \frac 1{z^2}(
{1-e^{-z \partial_x} })^2 - 2 m \partial_t\right)
\phi(x,t)=0 .
\ee
 By construction, the Galilei generators $\{K,H,P,M\}$ are
symmetries of $E_z$.   For the dilation we obtain from the  commutation
rules (\ref{bb}) that $ [E_z,D]= 2 E_z  $, 
so it is also a symmetry. For the conformal transformation we have
\be
 [E_z,C]=- K(1- e^{-2 z P})/z + M - 2 M D - 2 z MKH .
\ee 
By introducing the quantum realization  (\ref{bc}) we  get
\be 
 [E_z,C]=\left\{ t (e^{-z \partial_x} +1) - z m x e^{z \partial_x}\right\}
E_z + m (1+2a) .
\ee
Therefore  we have to set 
$a=-1/2$
 in order to obtain $C$ as a symmetry.

\section{A time discretization of the Schr\"odinger 
equation}

We consider now the quantum Schr\"odinger algebra \cite{twob}
$U_z^{(t)}({\cal S})$, 
 associated with the (nonstandard) classical
$r$-matrix $r= z H\wedge (D +\frac 12 M)\equiv z H\wedge D'$:
\be
\begin{array}{l}
 \Delta(H)=1\otimes H + H \otimes 1\qquad
 \Delta(M)=1\otimes M + M \otimes 1 \cr
 \Delta(P)=1\otimes P + P\otimes e^{-zH}\cr
 \Delta(K)=1\otimes K + K\otimes e^{zH}-z D' \otimes
e^{2zH}P\cr
 \Delta(D)=1\otimes D + D\otimes e^{2zH} +  \frac 12 M\otimes
(e^{2zH}-1) \cr
 \Delta(C)=1\otimes C + C\otimes  e^{2zH}-\frac z2 D'\otimes  e^{2zH}  M ,
\end{array} 
\label{ca}
\ee
\be
\begin{array}{l}
 [D,P]=-P \qquad   [D,K]=K\qquad
 [K,P]=M  \qquad  [M,\,\cdot\,]=0\cr
  [D,H]= \frac 1z({1-e^{2zH}}) \qquad   [D,C]=2 C + z  (D')^2\qquad 
[H,P]=0\cr
 [H,C]=  D+\frac 12 M (1-e^{2zH}) \qquad  [K,C]=\frac z2(DK+KD+KM)\cr
   [P,C]=-K-\frac z2  (DP+PD+PM) \qquad  [K,H]= e^{2zH}P .
\end{array} 
\label{cb}
\ee
 
A realization of $U_z^{(t)}({\cal S})$ 
reads ($b= m/2-2$):
\be
\begin{array}{l}
H=\partial_t \qquad  P=\partial_x \qquad M=m \cr
K=- (t + 2z) e^{2 z \partial_t}  \partial_x  - m x  \qquad  
 D= (t+ 2z) \frac 1z ({e^{2 z \partial_t}-1})+
x\partial_x - a   \cr
C=(t^2 - 2z b t) \frac 1{2z}({e^{2 z \partial_t}-1})  + t x
\partial_x - a t  + \frac m2 x^2    - 2 z (b +1) e^{2 z
\partial_t}  \cr
\qquad - \frac z2 x^2 \partial_{x}^2  -  z (b-a + \frac 12) x
\partial_x -\frac z2(b-a)^2 .
\end{array} 
\label{cc}
\ee

The deformed Casimir for the Galilei  sector is
$E_z =P^2-  M \frac 1z({1-e^{-2zH}})$; its  action
  on $\phi(x,t)$ through (\ref{cc}) leads to a time
discretization of the SE:  
\be
  \left(\partial_x^2 -   m \frac 1z
({1-e^{-2z\partial_t}})\right)
\phi(x,t)=0 .
\ee

As in Sec. 1 the generators $\{K,H,P,M,D\}$ are
easily shown to be symmetries of this SE. For the conformal
transformation, we find that
\bea
&&[E_z,C]=-(KP+PK + 2 M D e^{-2 z H}) + M(M+2) (1-e^{-2 z H})\cr
&&\qquad\qquad\quad - z(DP^2 + 2 PDP + P^2 D+2 P^2 M)/2.
\eea
Now adopting  the quantum realization (\ref{cc}):
\bea
&&[E_z,C]=2(t+ 2z -  z x \partial_x) E_z -   z \left\{ m + 2
(1-a)\right\}\partial_x^2\cr
&&\qquad\qquad\quad +m(1+ 2 a e^{-2 z \partial_t})
+ m (m+2 ) ( 1 - e^{-2 z \partial_t}) ,
\eea
and   setting  again
$a=-\frac 12$ gives us
$[E_z,C]=2\left\{ t+  \frac z2 (1- m - 2   x \partial_x)\right\} E_z $.

Finally, we recall that there exist other quantum
SE \cite{FVa,FVb,Dobrev} although without a known associated
Hopf structure.

\section*{Acknowledgments}

A.B. and F.J.H. have been partially supported by DGICYT (Project 
PB94--1115) from the Ministerio de Educaci\'on y Cultura  de Espa\~na
and by Junta de
Castilla y Le\'on (Projects CO1/396 and CO2/297).
P.P. acknowledges A.E.C.I. for the award of a fellowship.

%\section*{References}

\end{document}